\documentclass[journal]{IEEEtran}
\usepackage{cite}
\usepackage{amsmath}
\usepackage{amssymb}
\usepackage{relsize}
\usepackage{mathtools}
\DeclarePairedDelimiter{\abs}{\lvert}{\rvert}
\usepackage{graphicx}
\graphicspath{ {figures/} }
\usepackage{epstopdf}
\usepackage{subfigure}
\usepackage{cleveref}
\usepackage{caption}
\usepackage{multirow}
\usepackage{tabularx}
\usepackage{makecell}
\usepackage{multicol}
\usepackage{booktabs} 
\usepackage{rotating}
\usepackage{float}
\usepackage{enumitem}
\usepackage{tikz}
\usepackage{chngcntr}
\usetikzlibrary{shapes,arrows}
\tikzstyle{decision} = [diamond, draw, fill=blue!20, 
    text width=5em, text badly centered, node distance=3cm, inner sep=0pt]
\tikzstyle{block} = [rectangle, draw, fill=blue!20, 
    text width=5em, text centered, rounded corners, minimum height=4em]
\tikzstyle{line} = [draw,very thick, -latex']

\ifCLASSINFOpdf
\else
\fi
\begin{document}
%
\title{Computationally Efficient Day-Ahead OPF using
Post-Optimal Analysis with Renewable and Load
Uncertainties}
%
%
%

\author{Parikshit~Pareek,
        and~Ashu~Verma}
\maketitle

\begin{abstract}
This paper presents a method to handle renewable source and load uncertainties in  Dynamic Day-ahead Optimal Power Flow (DA-OPF) using post-optimal analysis of linear programming problem. The method does not require the uncertainty distribution information to handle it. A new Participation Factor (PF) to distribute changes caused by uncertainty has been developed based on the current optimal basis. The proposed PF takes care of all the constraints and ensures optimality with uncertain renewable generation and load using Sensitivity Analysis (SA) and Individual Tolerance Ranges (ITR) for individual and multiple simultaneous changes respectively. For quantification of confidence level, standard density distribution of solar power output and load is used. The test results on IEEE 30-Bus establishes the applicability of the proposed method for handling single and multiple bus uncertainties.    
\end{abstract}


%
\IEEEpeerreviewmaketitle

\section{Introduction}
\IEEEPARstart{T}{he} Renewable Sources (RESs) have become a need of the hour due to global warming issues associated with conventional sources of energy. This makes the large-scale integration of renewable sources an important task in order to extract the complete benefit of these sources. Together, the uncertain and intermittent nature of these possess difficulties in operation and control of the electric grid. Along with these, the system load has also been posing uncertainty challenges due to the introduction of new dynamic components. Traditionally, generator base points are calculated  over a snapshot of data through Day-ahead Economic Dispatch (DED). Further, participation factors are calculated as the means of balancing the load and generation when it moves away from the base-point~\cite{wood1996}. Intra-day static economic dispatch or look-ahead economic dispatch is also used to increase the utilization of RESs~\cite{xie2011look}. Day-ahead Optimal Power Flow (DA-OPF) works on the objective to operate the power system economically while keeping various control parameters under permissible limits. The applicability of DA-OPF or DED varies with the market structure. They provide generator set points to meet the forecasted demand in a vertically integrated market while objective will be of profit maximization in an unbundled system where bids will be built upon the results of DA-OPF. Nevertheless, the uncertainty in load and RESs disrupt the power balance and affect the bidding process as well.

Historically, these problems are solved using approaches which come under the category of deterministic optimization and can provide infeasible solutions even with a small perturbation in balancing equation parameters. A safe DA-OPF method is required in order to get solutions which remain at least feasible, if not optimal, under uncertainty. Recently, a number of works on security constrained unit commitment and economic dispatch models using Stochastic Optimization (SO)~\cite{li2014mean}\cite{papavasiliou2011reserve}\cite{wang2008security}\cite{ryan2013toward}, chance constraint variants of optimization~\cite{bienstock2014chance}\cite{liu2010economic}, interval programming~\cite{wu2012comparison} and Robust Optimization (RO)~\cite{bertsimas2013adaptive}\cite{ding2016adjustable}\cite{jabr2015robust} has been published.

The stochastic and chance constraint optimization approaches shares a common drawback in the requirement of accurate Probability Density Functions (PDFs) of uncertain variables. In practice, solar and wind uncertainties are modeled as uni-modal and bi-modal distributions\cite{wang2011spinning} while load uncertainty is as normal distribution\cite{bo2009impact}. Parameter calculation of these distributions will require statistical models over the past data and will introduce the errors. Additionally, SO accuracy is highly dependent on the number of scenarios generated and the requirement is quite high even for higher values of duality gap\cite{papavasiliou2015applying}. The scenario reduction techniques\cite{dupavcova2003scenario} can reduce computational complexity but mostly at the cost of economic benefits \cite{wu2012comparison}\cite{papavasiliou2015applying}. On the other hand, RO and interval programming requires less information about the uncertainty parameters. RO handles uncertainty in a deterministic manner using a set-based description of it. The solution satisfies all the constraints with uncertainty i.e. provides a feasible solution. The result obtained from RO can be conservative in comparison to the SO \cite{zhao2013unified}. Meanwhile, all of these methods go on calculating the PF by some means or other. Works on RO, like \cite{ding2016adjustable}, modifies the objective function via including generator participation factors in it and formulating the adjustable robust OPF. However, this makes problem to fall out of the category of linear programming (LP) and issues with computational complexity comes into the picture, even more, when measures to reduce conservatism of RO like budget of uncertainty\cite{bertsimas2004price} and dynamic uncertainty set\cite{lorca2015adaptive} are applied. Similar strategy for calculating best-fit participation factor is used to handle variability over the scheduling period in real time dispatch problem \cite{reddy2015real}. This approach increases the size of the problem.  

Here, we take leverage from a vast volume of work available under Post-optimal or Sensitivity Analysis (SA) of Linear Programming (LP) problem\cite{dantzig2006linear} and apply it to handle solar power and load uncertainties. It works on the basis of continuity property of the nominal optimal solution. Initially, the LP problem is solved using a set of data and then the sensitivity of the results is checked under the various perturbations. There are methods available to handle situations like single value perturbation, conventional sensitivity analysis\cite{dantzig2006linear}, percentage change in the complete or a section of right hand side (RHS) vector of the problem, Tolerance Range (TR)\cite{ravi1989tolerance},Individual Tolerance Range (ITR)\cite{wondolowski1991generalization} and random changes in parameters, perturbation analysis\cite{arsham1990perturbation}. All the calculations are based on the validity of current optimal basis under perturbations. As the basis information is available with the solution itself,  uncertainty handling using this method has an essential advantage over the SO and RO in terms of computational complexity. This approach does not require any statistical information about the uncertainty and can be implemented in existing market environment easily. The issue of variation of Locational Marginal Price (LMP) due to uncertainty\cite{bo2009impact}, is also eliminated as the dual variables remain unchanged through this approach. A new Participation Factor (PF) is calculated for each of the generators to share the change in power occurred due to the uncertainty. The new solution will be an optimal solution not just feasible.  
This will help in incorporating the uncertain sources in the electricity market.\

Contributions of the present work can be summarized as: 
\begin{itemize}
\item Development of a method which ensures optimality under uncertainty and keeps the structure of the original DCOPF problem unchanged. This makes the method computationally efficient over present methods.  
\item This work proposes a new participation factor (PF) which can be used to distribute the change caused by uncertainty among generators. It includes all the constraints in providing the change in magnitude for each generator. 
\item Proposed approach does not require the probability distribution of uncertainty for handling it. Instead,  PDFs are used to quantify the confidence level in uncertainty handling.
\item The method can be applied to handle both load and RES uncertainty. It can handle multiple uncertainties over the day-ahead scheduled operating solutions as well. 
\end{itemize} 
\ \ Rest of the paper is organized as follows. Section II describes the DCOPF formulation and mapped Linear Programming (LP) problem. Section III is a detailed account of proposed PF and uncertainty handling approach, range update method and proof of LMP invariability. Section IV presents the test results of IEEE 30 Bus system for the proposed method and section V contains the conclusion. Details of solar and load uncertainty modeling are given in the appendix.   
\section{Problem Formulation}

\subsection{Dynamic DA-OPF Formulation}
The AC Optimal Power Flow (ACOPF) is a nonlinear, and nonconvex problem. Even for systems of small size, the ACOPF become a large scale problem which takes more time to solve than the current practices of day-ahead and real time market allow~\cite{cain2012history}. Amongst various approximate formulations, DC Optimal Power Flow (DCOPF) or DC Economic Dispatch (DCED) is the one used by most of the Independent System Operators (ISOs) across the globe. Current solvers provide an efficient result with the admittance matrix 
model of DCOPF\cite{stott2009dc} and negates the problem of loss factor based DCOPF implementation\cite{litvinov2004marginal} related to slack bus selection. A recent work \cite{brent2017} provides a clear account of the current practices related to DCOPF in market operations. It also reiterates the importance of current linear formulation based approaches in the context of their advantage, over other complex formulations, in terms of computational performance and effective integration with economic theory~\cite{caramanis1982optimal}. The formulation for generation cost minimization objective with constraints over line flow limits and power balance is given as: 
\begin{subequations}\label{eq:dcopf}
\begin{gather}
minimize \mathlarger{\sum_{k\in \Omega_g }} C_f(P_{g{_k}}) \\
subject\ to \mathlarger{\sum_{k\in \Omega_{G{_k}}}} P_{g{_k}}-P_{l_{i}}=\mathlarger{\sum_{j\in \Omega_{ij}}} \dfrac{\theta_i-\theta_j}{X_{ij}}, \forall ~ i \in \Omega_j \label{constraints1}\\
-\overline{P}_{ij} \leqslant P_{ij}=\dfrac{\theta_i-\theta_j}{X_{ij}} \leqslant \overline{P}_{ij},\ \forall ~ ij\in \Omega_L \label{constraint2}  
\end{gather}
\end{subequations}\ \
Ramp rate constraints are taken to formulate the dynamic dispatch problem which in turn keeps turbine thermal gradients within the limit. They are modeled with the knowledge of past hour generation as ($\forall \ k\in \Omega_G$):
\begin{equation}\label{constraint3}
max [\underline{P}_{g{_k}},P^{h-1}_{g{_k}}-R^{down}_{G{_k}}] \leqslant P_{g{_k}} \leqslant min [\overline{P}_{g{_k}},P^{h-1}_{g{_k}}+R^{up}_{G{_k}}]
\end{equation}
\ \ Here, $C_f(P_{g{_k}})$ is linear cost and ${P}_{g{_k}}$ is generation of $k^{th}$ generator, ${P}_{l{_i}}$ is load at $i^{th}$  bus, $X_{ij}$ is reactance of the line between buses $i$ and $j$, $\overline{P}_{ij}$ is capacity of the line between buses  $i$ and $j$, $\theta_i$ is the angle at bus $i$, $\Omega_G$ is generator index, $\Omega_L$ is line index, $\Omega_j$ is bus index, $\Omega_{G{_k}}$ is generators at bus $k$, $R^{up}_{G{_k}}$ and $R^{down}_{G{_k}}$ are ramp rates in both directions, $\overline{P}_{g{_k}}$ and $\underline{P}_{g{_k}}$ is maximum and minimum generation capacity of generator $k$ respectively. 
\subsection{LP Problem Formulation}
The Linear Programming Problem (LPP) to solve the Dynamic DA-OPF is modeled by converting all the constraints as equality constraints to facilitate the post-optimal analysis. $P_{G_{k}}^s$ are control variables and the slack bus is not selected beforehand. The LP problem looks like: 
\begin{subequations}\label{eq:LPP}
	\begin{gather}
	minmize\ Z=\boldsymbol{C'X}\\
	subject\ to\ \boldsymbol{AX=b}\\
\quad{\qquad}[\boldsymbol{\theta  S}]\geqslant \boldsymbol 0
	\end{gather}
\end{subequations}
Here, $\boldsymbol{X=[P_{g_{k}} \theta S]'}$, $\boldsymbol{b=[\overline{P}_{g{_k}} -\underline{P}_{g{_k}} \overline{P}_{ij} -\overline{P}_{ij} P_l]'}$. Constraint matrix $\boldsymbol{A}$ is made up of all the constraints (\ref{constraints1}, \ref{constraint2} and \ref{constraint3}) and sparse in nature. The solution of this problem will give $\boldsymbol{B}$ as optimal basis matrix with $i_b$ as basis index, $\boldsymbol{X_b=B^{-1}b}$ the basic optimal solution.  Bold symbols represent matrices and vectors wherever written. The uncertainty handling method proposed in this paper is based upon the continuity of $\boldsymbol{B}$ over the perturbations in $\boldsymbol{b}$. Bus angle value ($\theta_i, \forall i \in \Omega_j$) is calculated by subtracting slack bus angle from $\boldsymbol{\theta}$ obtained in solution of LP (\ref{eq:LPP}).

	


%
\section{Proposed Method for Uncertainty Handling}\label{sec:proposed}
Conventionally, for matching the load and generation over day-ahead scheduled optimal solution, the participation factors of each generator is calculated based upon the concept of moving the generators from an economic optimal to another one to share the load change \cite{wood1996}. However in the work, the proposed participation factor, $\beta_k$, is developed based upon the current optimal basis obtained with the day-ahead solution for handling the uncertainties in solar PV output and load without changing the problem structure. Mathematically, 
 \begin{equation}\label{beta}
 \beta{_k}=\dfrac{\Delta P_g{_k}}{\sum_{j=1}^{nbus}\Delta P_l{_j}}=\frac{\sum_{j=1}^{nbus}R(k,j)~ \Delta P_l{_j}}{\sum_{j=1}^{nbus}\Delta P_l{_j}},\\
 \end{equation}
 \begin{subequations}\label{eq:PF}
 	\begin{gather}
 	\Delta P_g{_k} \in \big \{ \Delta X_b(1),\Delta X_b(2),......\Delta X_b(ng) \big \}\\
 where,\  \Delta X_b(k)=\sum_{i=1}^{m}B^{-1}(k,i)~ \Delta b(i),\\
 or\quad   \Delta X_b(k)=\sum_{i=1}^{nbus}R(k,i)~ \Delta P_l{_i},
 	\end{gather}
 \end{subequations}
 \ \ Here, $m$ is number of constraints, $\boldsymbol{R}^{m\times nbus}$  matrix having columns of $\boldsymbol{B}^{-1}$ corresponding to the load values. The $\beta_k$ is participation factor ($ \forall k \in \Omega_G$). It is the ratio of change in $k^{th}$ generator schedule to the total change in load or PV output. The $\Delta P_g{_k}$ quantifies the change in the generator output $P_g{_k}$ ( $ \forall k \in \Omega_G$). DCOPF (\ref{eq:dcopf},\ref{constraint3}) is mapped in LP problem (\ref{eq:LPP}) such that first $ng$ variables of $X_b$ will give the generator set points and changes in them (\ref{eq:PF}a). To obtain $\Delta P_g{_k}$, continuity property of $\boldsymbol{B}$ is used according to which $\beta_k$ is valid under that changes in the day-ahead problem parameters till $\boldsymbol{B}$ maintains its continuum.\
  
 The single value perturbations like PV output, in DCOPF (\ref{eq:dcopf}) parameters, conventional sensitivity analysis (SA) is used\cite{dantzig2006linear}. The feasibility and optimality conditions of LP solution are:
 \begin{equation}\label{feasi}
 \boldmath{B^{-1}b} \geqslant\boldmath{0}
 \end{equation}
 \begin{equation}\label{eq:opti}
 \boldmath{C'-C'_BB^{-1}A} \geqslant\boldmath{0}
 \end{equation}
 \ \ As the PV output can be modeled as negative load, the uncertainty in PV output is treated as changes in the RHS vector $\textbf{\textit{b}}$ which only appears in the feasibility condition (\ref{feasi}). Therefore, the uncertainty must satisfy only feasibility condition to remain optimal as well. Now, with perturbations in $i^{th}$ element value and $\boldsymbol{e_i}$ being a $i^{th}$ column of $\boldsymbol{B}^{-1}$:
 \begin{equation}
 \boldmath{B^{-1}(b+\Delta b e_i)} \geqslant\boldmath{0}
 \end{equation}
\ \ With $\boldsymbol{f} $ as  $i^{th}$ column of $\boldsymbol B^{-1}$ and $\alpha_{ij}$ as $j^{th}$ element of the $\boldsymbol f $ then:
\begin{equation}
\boldmath{X_b+\Delta bf} \geqslant\boldmath{0}
\end{equation}
\begin{equation}\label{currentoptsol}
X_{b}(j)+\boldsymbol{\Delta b}\alpha_{ij} \geqslant 0
\end{equation}

Unvaryingly, 
\begin{subequations}\label{deltab}
\begin{gather}
\begin{split}
\bigg\{ max ({-X_{b}(j)}/{\alpha_{ij}})for ~ j|\alpha_{ij}>0 \bigg \} \geqslant \mathlarger{\Delta b(j)} \\ \geqslant \bigg \{ min({-X_{b}(j)}/{\alpha_{ij}}) for~ j|\alpha_{ij}<0\bigg \}
\end{split}\\
\begin{split}
\bigg \{ min({-X_{b}(j)}/{\alpha_{ij}}) for~ j|\alpha_{ij}<0\bigg \} \geqslant \mathlarger{\Delta P_{PV}(j)} \\ \geqslant \bigg\{ max ({-X_{b}(j)}/{\alpha_{ij}})for ~ j|\alpha_{ij}>0 \bigg \} 
\end{split}\\
P_{PV}(j)+\Delta P_{PV}^{min}(j) \geqslant P_{PV}^{new}(j) \geqslant P_{PV}(j)+\Delta P_{PV}^{max}(j)
\end{gather}
\end{subequations}
\ \ Here, $P_{PV}(j)$ and  $\Delta P_{PV}(j)$ are present value and change in $P_{PV}$ installed at $j^{th}$ bus of the system. $\Delta P_{PV}^{min}(j)$ and $\Delta P_{PV}^{max}(j)$ are the minimum and maximum change up to which the $\boldsymbol{B}$ is valid as per the feasibility condition (\ref{feasi}). The calculations of limits are based upon the ratio of negative of current optimal solution to the $\alpha_{ij}$ element of $\boldsymbol{B}^{-1}$. Minimum and maximum values of this ratio with condition on $j$ ($\alpha_{ij}< 0$ and $\alpha_{ij}>0$)are taken to get the boundary values. \

In case load/generation uncertainty needs to be handled simultaneously at more then one bus ITR \cite{wondolowski1991generalization} is used which provides features like 1) simultaneous and independent load modifications are permitted; 2) easy calculation from optimal solution information; 3)even valid in the occurrence of primal or dual degeneracy; 4) dual variables does not vary after handling the uncertainty. Concentrating upon the load change in $j^{th}$ bus as fraction, $\delta_{P_l}(j)$, is given as:
\begin{equation}\label{delp}
  \delta_{P_l}(j)=\dfrac{X_b(j)}{\mathlarger{\sum_{k\in \Omega_j }}{\abs {\alpha_{jk}\boldsymbol P_l(k)}}}~\forall j\in \Omega_j
\end{equation}
\ \ Further, the permissible load increase ($\Delta P_l^+(j)$) and decrease ($\Delta P_l^-(j)$) at $j^{th}$ bus is given by minimum values of multiplication of absolute of load ($P_l(j)$) and $\delta_{P_l}(k)$ with conditions on $k$ as:
 \begin{equation}
 \Delta P_l^+(j) = min\big\{(\abs{P_l(j)}\times\delta_{P_l}(k))\mid  \alpha_{jk}<0\}
 \end{equation}
 \begin{equation}
 \Delta P_l^-(j)  = min\big\{(\abs{P_l(j)}\times\delta_{P_l}(k))\mid  \alpha_{jk}>0\}
 \end{equation}
 \ \ Therefore, the new load $P_l^{new}(j)$ must follow the relation ($\forall ~j\in \Omega_j$):
 \begin{equation}\label{eq:itr}
 P_l(j)+\Delta P_l^-(j)\leqslant P_l^{new}(j)\leqslant P_l(j)+\Delta P_l^+(j)
 \end{equation}
\ \ Here at situations of primal or dual degeneracy the value of 
$\delta_{P_l}(j)$ (\ref{delp}) is set to zero as any perturbation will lead to change in the optimal basis \cite{wondolowski1991generalization}. Thus, the  $\beta{_k}$ can be obtained for both $\Delta P_{PV}$ and $\Delta P_l$ from (\ref{beta}).  



In practice, due to uncertainties in solar PV generation and loads, deviations from day-ahead operating solution in a system should be addressed. Further, with an increase in RESs, a system can have multiple solar PV systems. To handle these cases, the SA and ITR values are proposed to be modified successively. Conceptually, after handling each uncertainty, the day-ahead optimal solution ($\boldsymbol{X_b^{0}}$) moves to ($\boldsymbol{X_b^{1}}$) the condition (\ref{currentoptsol}) will become:
\begin{subequations}
	\begin{gather}
X_{b}^{1}(j)+\boldsymbol{\Delta b^1}\alpha_{ij} \geqslant 0\\
where,\ \boldsymbol{X_{b}^{1}=B^{-1}(b^o+\Delta b^o)}  
	\end{gather}
\end{subequations} 

\ \ Hence, updated ranges are obtained (\ref{deltab}) and (\ref{delp}-\ref{eq:itr}) to give the new ranges of $\Delta P_{PV}$ amd $\Delta P_l$  up to which uncertainties are further allowed under current optimal basis. This range update method makes proposed approach apt for handling the uncertainties in the load and RESs over a scheduling period in real-time. In order to showcase the applicability of the proposed method, the confidence level is calculated as area under the PDFs within the limits (\ref{deltab},\ref{eq:itr}) and it quantifies the degree of robustness simultaneous.  

\subsection{LMP Invariability}\label{sec:LMP} 
The LMP for the lossless DCOPF formulation is made-up of components of energy and transmission congestion and given as\cite{li2007dcopf} (~$ \forall j\in \Omega_j$):
 \begin{subequations}
\begin{gather}
LMP_j=LMP^{energy}+LMP^{congestion}_j \\  
 LMP^{energy}=\lambda \\
 LMP^{congestion}_j=\sum\limits_{k \in \Omega_{L}}GSF_{k-j}\times \mu_k
\end{gather}
  \end{subequations}
\ \ Where, $\lambda$ is the dual variable corresponding to the power balance constraint (\ref{eq:dcopf}b), $\mu_k$ is the dual variable corresponding to the thermal constraints of each line (\ref{eq:dcopf}c) and $GSF_{k-j}$ is the generator shift factor to line $k$ from bus $j$. As the $GSF$ depends upon the network\cite{wood1996}, it can be taken as constant over the time period under present study. Therefore, the LMP depends upon the dual variables obtained from the solution of the LP problem (\ref{eq:LPP}) corresponding to the Dynamic DA-OPF (\ref{eq:dcopf},\ref{constraint3}). 
The dual optimal solution corresponding to the primal one is given as:
\begin{equation}\label{eq:dual}
\boldsymbol{\pi=C'_bB^{-1}}
\end{equation}\
As the SA and ITR work on the invariability of the current optimal basis, both $\boldsymbol C_b$ and $\boldsymbol B^{-1}$ remains constant within the uncertainty range which can be handled by the proposed method. Therefore, both $\lambda$ and $\mu_k$ remain constant keeping the LMP value invariable.       

\section{Results and Discussion}
\begin{figure}[b]
 	\includegraphics[width=\linewidth,height=4.5cm]{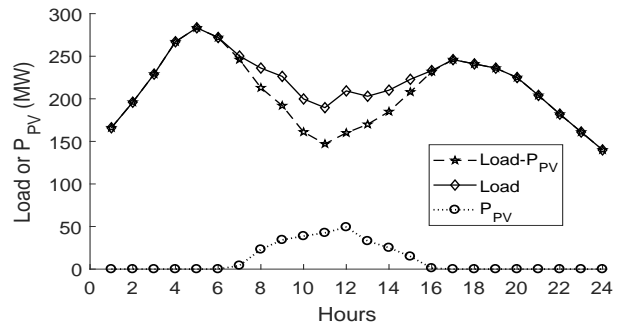}
 	\caption{Load profiles used for the system}
 	\label{fig:load}
\end{figure}

The IEEE 30 Bus test system \cite{alsac1974optimal} is used to show the applicability of the proposed method for handling uncertainties in solar PV generation and load. It has 41-branches and 6-generators at Bus 1, 2, 5, 8, 11, and 13. Solar PV is taken of 60 MW peak and placed at different buses to show the variable robustness of network against uncertainties. The load profiles and the PV output is shown in Fig.\ref{fig:load}. The solar PV output starts at $7^{th}$ hour and lasts up to $16^{th}$ hour. The peak is occurring at $12^{th}$ hour as 49.3 MW. For load, a peak is 283.4 MW occurring at $4^{th}$ hour. The dynamic DA-OPF interval is taken as 1 hour which can easily be taken as 15 min. or so for real-time economic dispatch problem. MATLAB R2017a on an Intel i5 7200U processor system with 8 GB memory is used for the tests.         
\

\begin{figure}[b]
	\centering
 	\includegraphics[width=\linewidth,height=4.5cm]{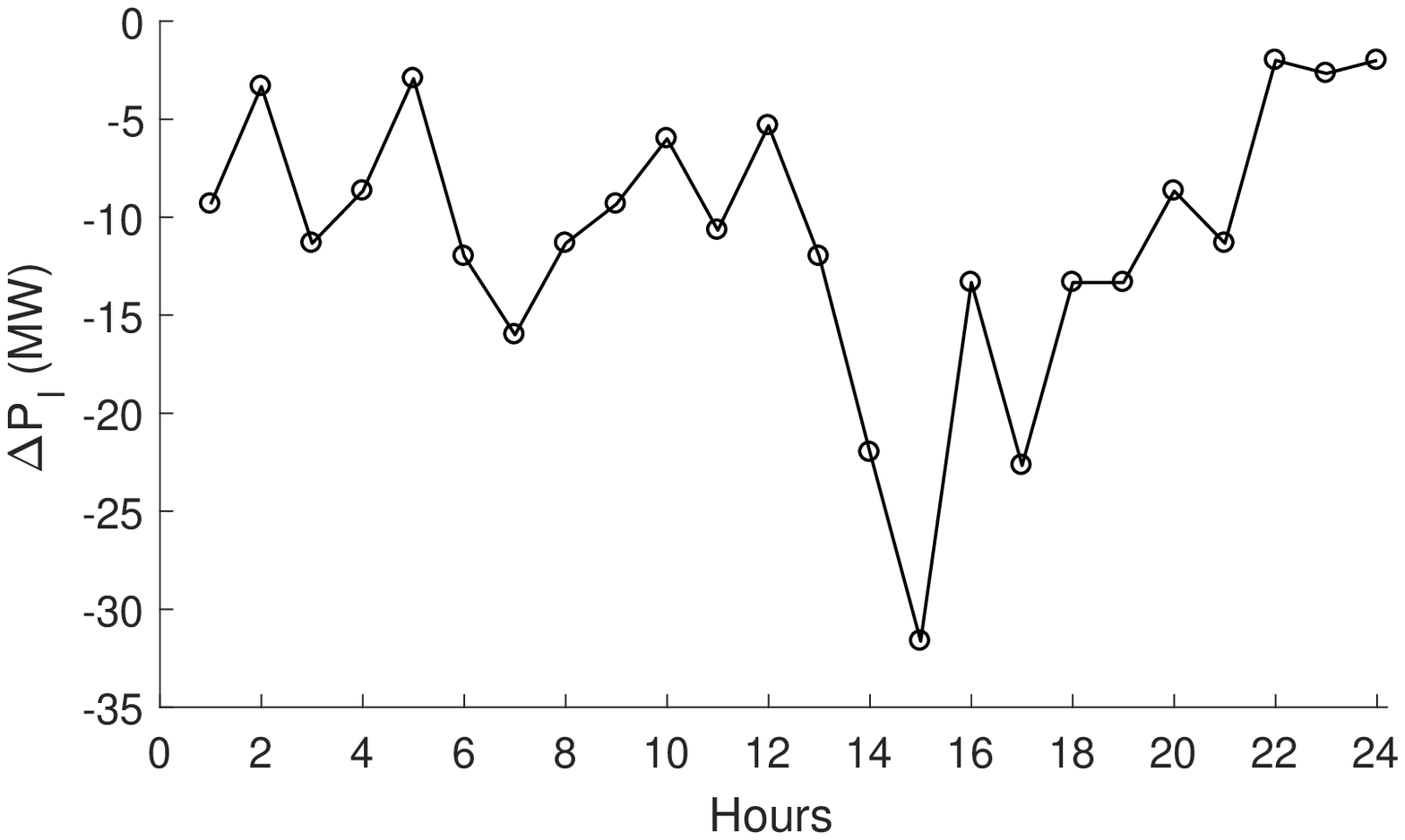}
 	\captionsetup{justification=centering}
 	\caption{Uncertainty caused system load change}
 	\label{fig:SF_load}
\end{figure}

The load data (Fig.\ref{fig:load}) used for DA-OPF solution can get changed due to uncertainty. One such situation is shown in Fig.\ref{fig:SF_load} where the load decreases due to uncertainty. This change in load is distributed among buses according to their load factors.Now the $\beta_{k} $, as defined in section \ref{sec:proposed},  gives the factor by which different generators will share the changed load. Table \ref{tb:SF} shows the power generated by each generator before and after the load uncertainty.\

Results indicate that the two most costly generators (generator no. 4 and 5) never shares any decreased load as they are already at minimum capacity. Generator no. 3 and 6 shared the load at two instances each and at a majority of the instances, the load is shared by the first and second generator. The $\beta_k$ does not only considers economy but considers ramp-rate constraints as well in the same way they are taken in the DED problem. Generation and transmission limits are also respected and hence the results obtained are similar to the onces obtained by re-optimization with changed loads. 

\subsection{Solar PV Uncertainty Handling}
\begin{figure}[t]
 	\centering
  	\includegraphics[width=\linewidth,height=4.5cm]{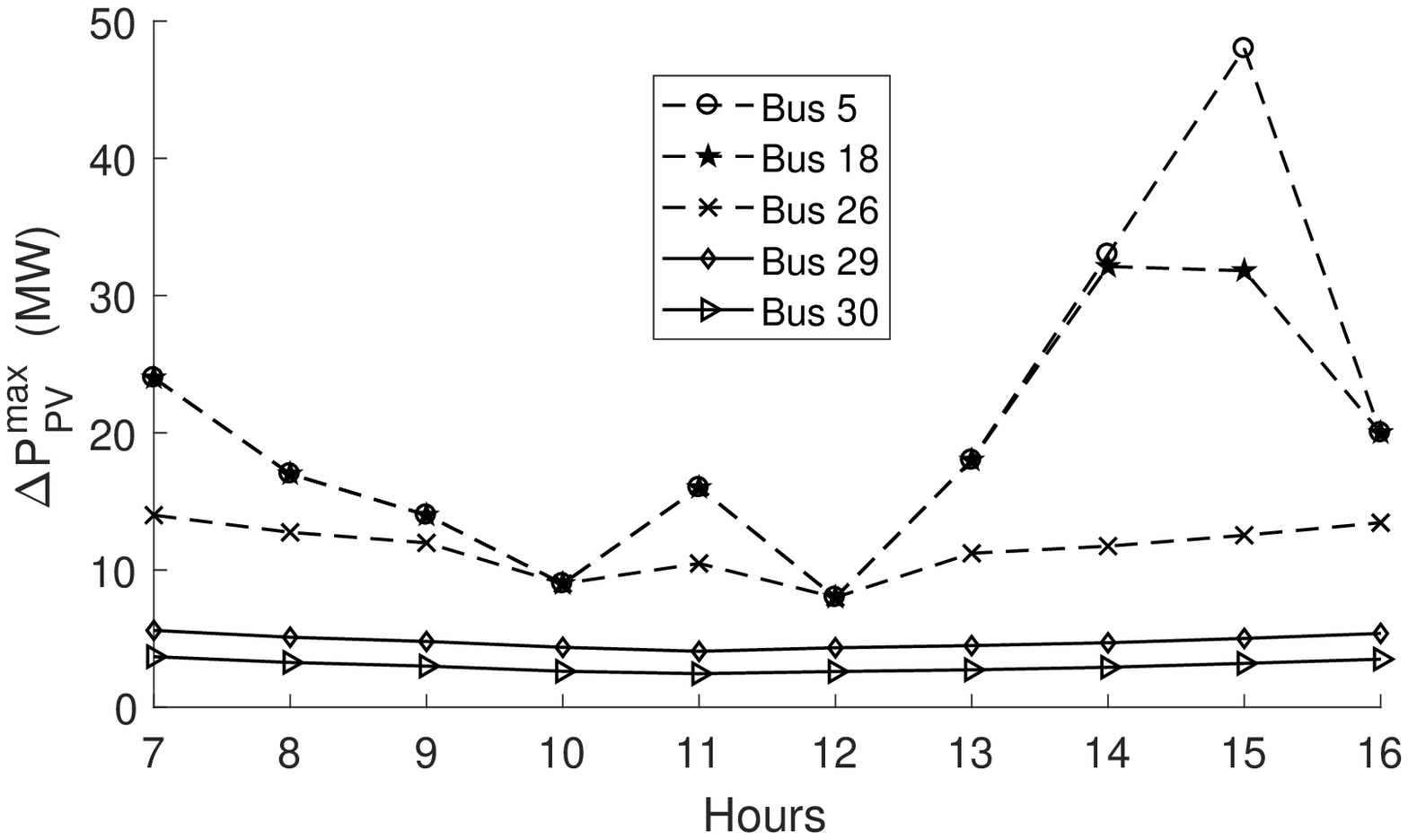}
  	\captionsetup{justification=centering}
   	\caption{Acceptable increase PV output due to uncertainty}
  	\label{fig:negi_SA}
 \end{figure}

\begin{figure}[t]
 	\centering
  	\includegraphics[width=\linewidth,height=4.5cm]{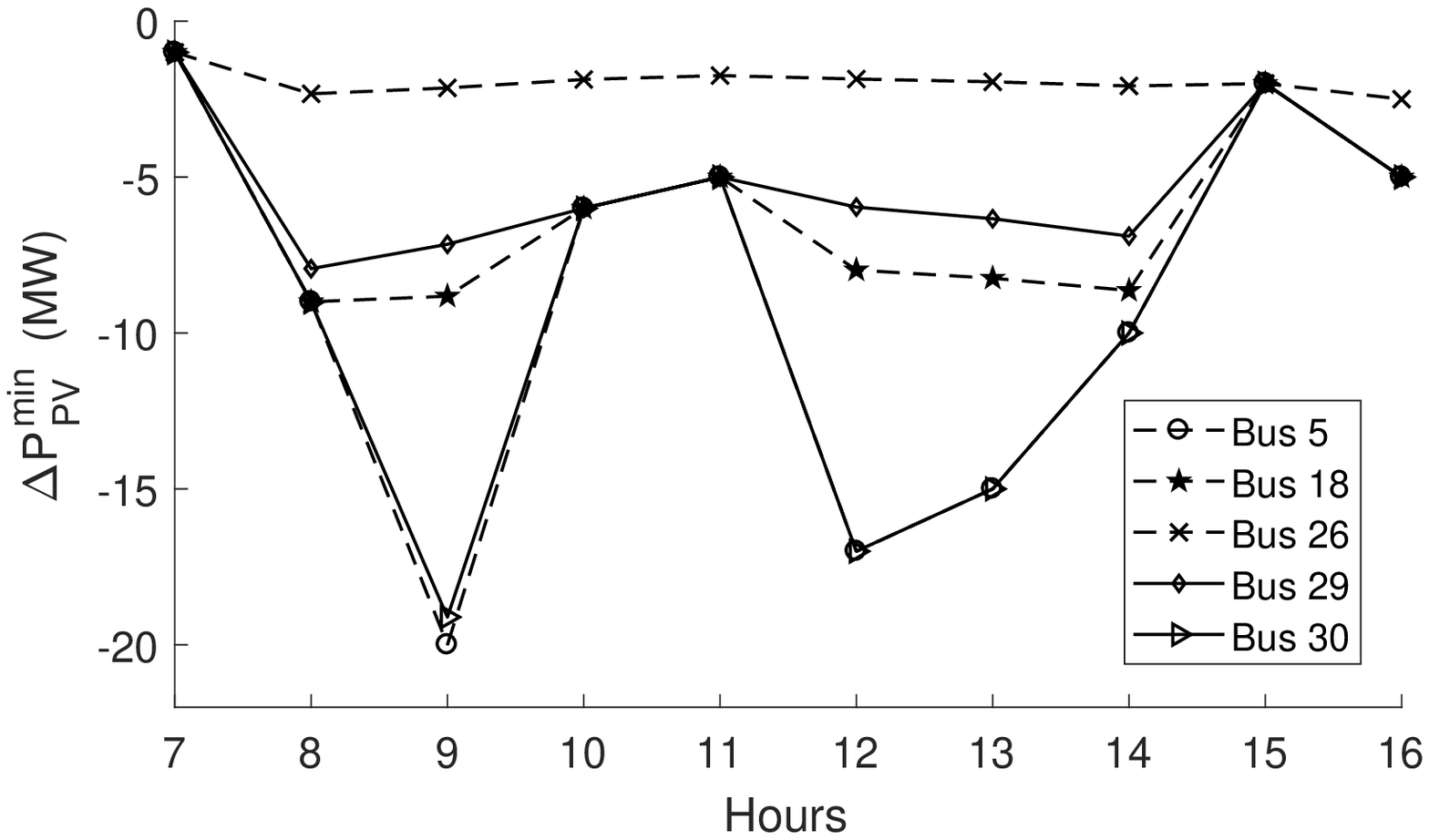}
  	\caption{Acceptable decrease PV output due to uncertainty}
  	\label{fig:posi_SA}
 \end{figure}

\begin{figure}[t]
	\centering
	 \begin{minipage}[b]{\linewidth}	 	
	 \centering
	  	\includegraphics[width=\linewidth,height=3.2cm]{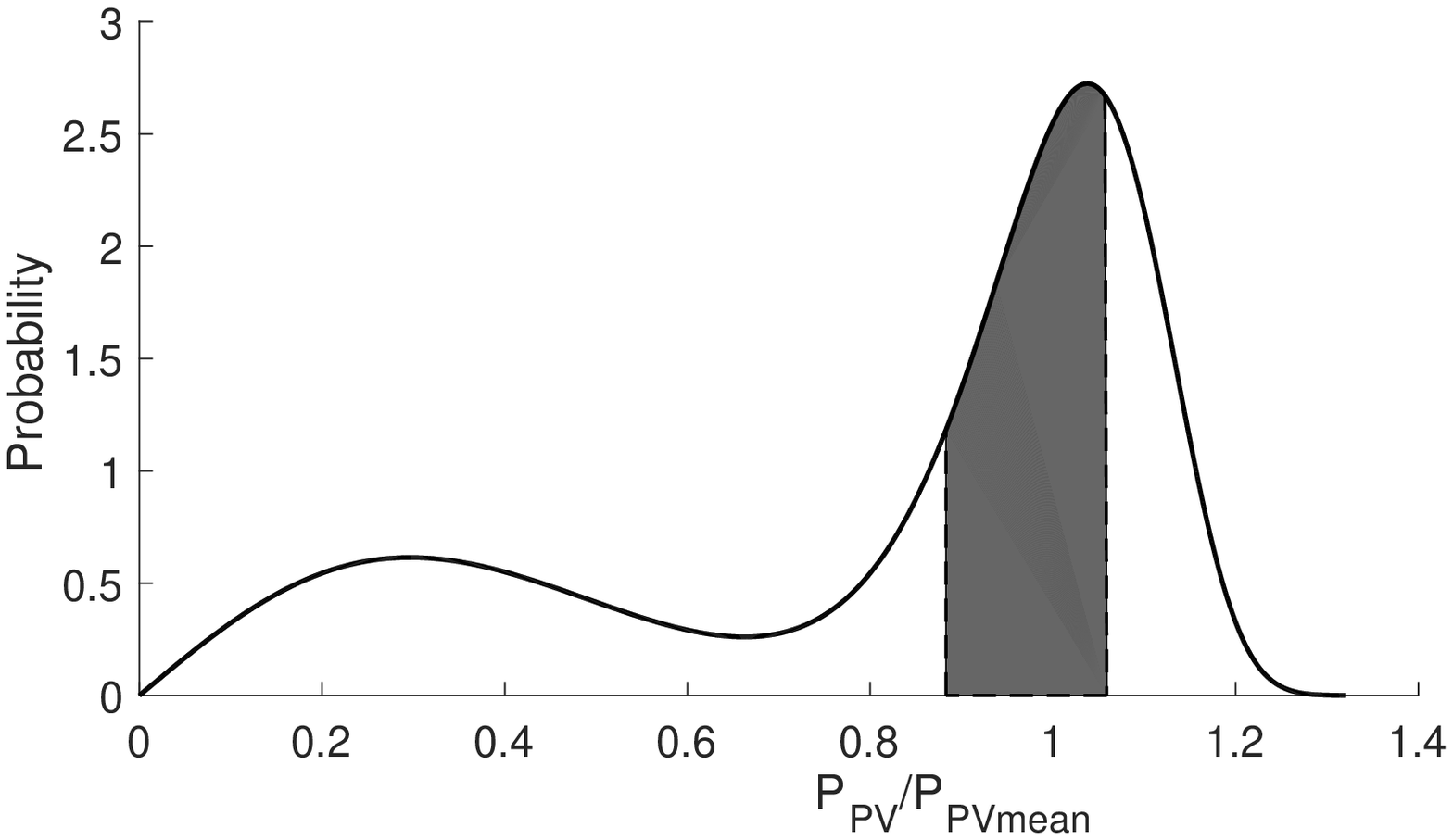}
	  	\\ \scriptsize{(a)}
     \end{minipage}
 \begin{minipage}[b]{\linewidth}
 	\centering
 	 	\includegraphics[width=\linewidth,height=3.2cm]{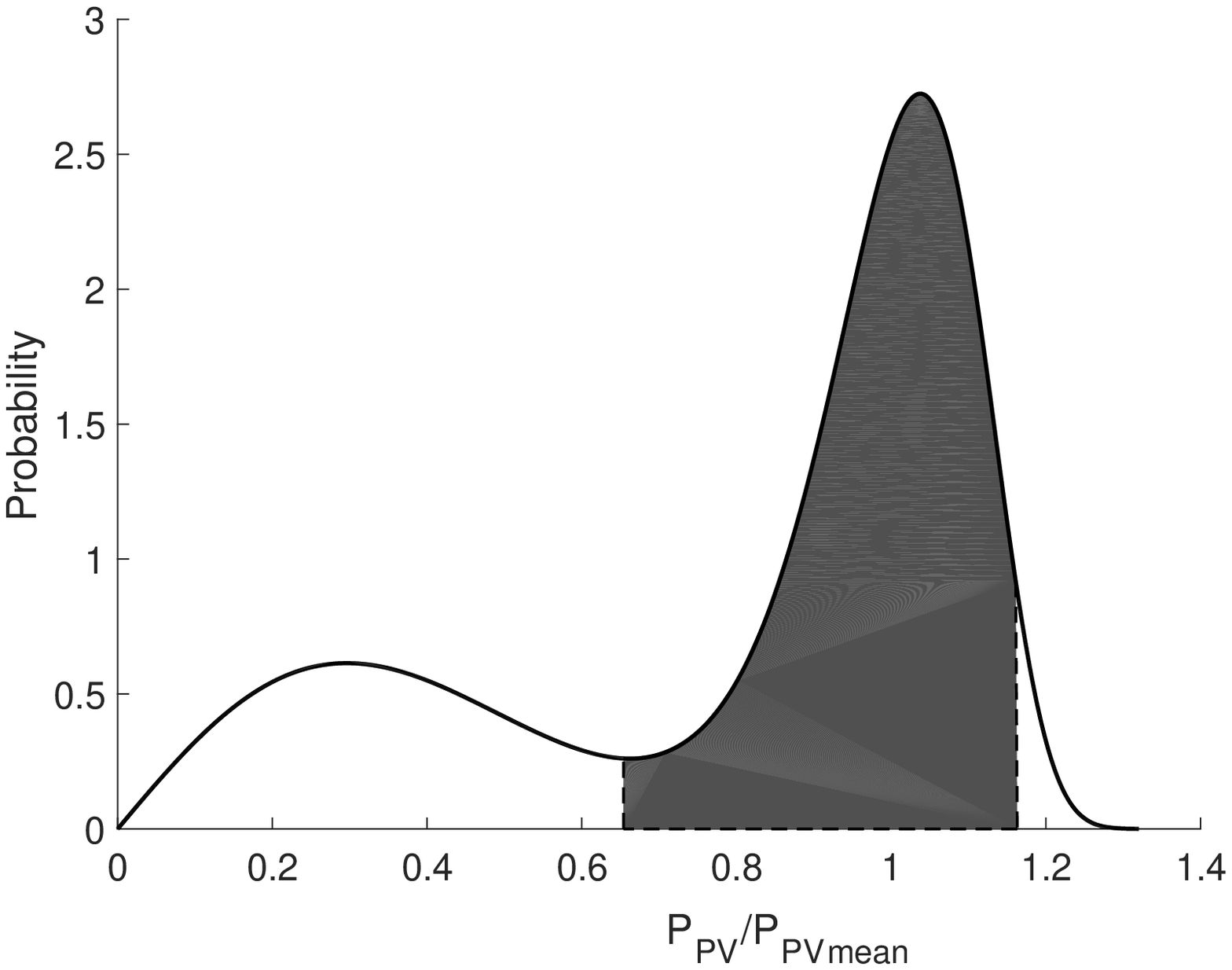}
 	\\ 
 	\scriptsize{(b)}
 	 \end{minipage}
  \caption{Area under the PDF (a) Bus no. 30 $11^{th}$ hour (b) Bus no. 5 $12^{th}$ hour }
  \label{fig:areafill}
\end{figure}

\begin{figure}[t]
	\centering
 	\includegraphics[width=\linewidth,height=4.5cm]{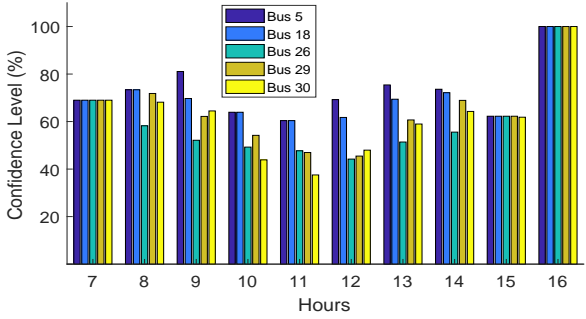}
 	\caption{Confidence level in PV uncertainty handling}
 	\label{fig:pv_area}
\end{figure}

The PV output uncertainty is handled by carrying out by single perturbation at the bus where PV is connected, as described in the section \ref{sec:proposed}, The acceptable values of $\Delta P_{PV}^{max}(j)$ and $\Delta P_{PV}^{min}(j)$ are shown in Fig.\ref{fig:negi_SA} and Fig.\ref{fig:posi_SA}. In order to quantify the confidence level, with respect to Fig. \ref{fig:un}(a), the area under the PDF curve within the acceptable $\Delta P_{PV}$ (\ref{deltab}b) is calculated as shown in Fig.\ref{fig:areafill}. The shaded dark strip is the area of the PDF which falls under the limits obtained through SA. Fig.\ref{fig:pv_area} depict the dependency of confidence level over the location of the PV plant. 
It can be observed from Fig.\ref{fig:pv_area} that confidence level is less at peak generation hours as the variation will be high in magnitude (MW) at that time. It proves that present method works on absolute MW limits and at peak generation uncertainty will be higher and confidence level will be lover. 
 The limiting factor for uncertainty handling at different buses is the power evacuation capacity and current load. The total evacuation capacity is 200 MW at $5^{th}$ bus while it is 16 MW for the $29^{th}$ bus. Similarly, $5^{th}$ bus has highest load fraction as well.It is clear from the results that confidence level is higher at $5^{th}$ Location of generators, load and evacuation capacity at adjoining buses is also influencing these values. Thus, the results present holistic nature of the proposed method for handling uncertainties.    
\def\arraystretch{1.1}  
\begin{table*}[t]
\caption{Generations before and after the uncertainty handling}
\label{tb:SF}
\resizebox{\linewidth}{!}{%
\begin{tabular}{c|cccccc|cccccc}
\hline \textbf{Hours} & \multicolumn{6}{c}{\textbf{Before Uncertainty}} & \multicolumn{6}{c}{\textbf{After Uncertainty}}                                \\

                & $P_{g_1}$ & $P_{g_2}$ & $P_{g_3}$ & $P_{g_4}$ & $P_{g_5}$ & $P_{g_6}$ & $P_{g_1}$ & $P_{g_2}$ & $P_{g_3}$ & $P_{g_4}$ & $P_{g_5}$ & $P_{g_6}$           \\
               \hline 
1              & 50.00    & 34.00    & 50.00    & 10.00    & 10.00    & 12.00      & 50.00              & \textbf{24.66} & 50.00              & 10.00   & 10.00   & 12.00              \\
2              & 55.00     & 59.00     & 50.00     & 10.00     & 10.00     & 12.00       & \textbf{51.66} & 59.00              & 50.00             & 10.00   & 10.00   & 12.00              \\
3              & 67.00     & 80.00     & 50.00    & 10.00    & 10.00    & 12.00      & \textbf{55.66} & 80.00             & 50.00             & 10.00   & 10.00   & 12.00              \\
4              & 92.00     & 80.00     & 50.00    & 10.00    & 10.00    & 25.00      & 92.00             & 80.00             & 50.00             & 10.00   & 10.00   & \textbf{16.33} \\
5              & 117.00    & 80.00     & 50.00    & 10.00    & 10.00    & 16.40   & 117.00            & 80.00             & 50.00             & 10.00   & 10.00   & \textbf{13.46} \\
6              & 110.00    & 80.00     & 50.00    & 10.00    & 10.00    & 12.00      & \textbf{98.00}    & 80.00             & 50.00             & 10.00   & 10.00   & 12.00              \\
7              & 85.00     & 79.00     & 50.00    & 10.00    & 10.00    & 12.00      & 85.00             & \textbf{63.00} & 50.00             & 10.00   & 10.00   & 12.00              \\
8              & 60.00     & 71.00     & 50.00    & 10.00    & 10.00    & 12.00      & 60.00             & \textbf{59.66} & 50.00             & 10.00   & 10.00   & 12.00              \\
9              & 50.00     & 60.00     & 50.00    & 10.00    & 10.00    & 12.00      & 50.00             & \textbf{50.66} & 50.00             & 10.00   & 10.00   & 12.00              \\
10             & 50.00     & 35.00     & 44.00    & 10.00    & 10.00    & 12.00      & 50.00             & 35.00             & \textbf{38.00}    & 10.00   & 10.00   & 12.00              \\
11             & 50.00     & 20.00     & 45.00    & 10.00    & 10.00    & 12.00      & 50.00             & 20.00             & \textbf{34.33} & 10.00   & 10.00   & 12.00             \\
12             & 50.00     & 28.00     & 50.00    & 10.00    & 10.00    & 12.00      & 50.00             & \textbf{22.66} & 50.00             & 10.00   & 10.00   & 12.00              \\
13             & 50.00     & 38.00    & 50.00    & 10.00    & 10.00    & 12.00      & 50.00             & \textbf{26.00}    & 50.00             & 10.00   & 10.00   & 12.00              \\
14             & 50.00     & 53.00     & 50.00    & 10.00    & 10.00    & 12.00      & 50.00            & \textbf{31.00} & 50.00             & 10.00   & 10.00   & 12.00              \\
15             & 50.00     & 76.00     & 50.00    & 10.00    & 10.00    & 12.00      & 50.00             & \textbf{44.36} & 50.00             & 10.00   & 10.00   & 12.00              \\
16             & 70.00     & 80.00     & 50.00    & 10.00    & 10.00    & 12.00      & \textbf{56.66} & 80.00             & 50.00             & 10.00   & 10.00   & 12.00              \\
17             & 84.00     & 80.00     & 50.00    & 10.00    & 10.00    & 12.00      & \textbf{61.33} & 80.00             & 50.00             & 10.00   & 10.00   & 12.00              \\
18             & 79.00     & 80.00     & 50.00    & 10.00    & 10.00    & 12.00      & \textbf{65.66} & 80.00             & 50.00             & 10.00   & 10.00   & 12.00              \\
19             & 74.00     & 80.00     & 50.00    & 10.00    & 10.00    & 12.00      & \textbf{60.66} & 80.00             & 50.00             & 10.00   & 10.00   & 12.00              \\
20             & 63.00     & 80.00     & 50.00   & 10.00    & 10.00    & 12.00      & \textbf{54.33} & 80.00             & 50.00             & 10.00   & 10.00   & 12.00              \\
21             & 50.00     & 72.00     & 50.00    & 10.00    & 10.00    & 12.00      & 50.00             & \textbf{60.66} & 50.00             & 10.00   & 10.00   & 12.00              \\
22             & 50.00     & 50.00     & 50.00    & 10.00    & 10.00    & 12.00      & 50.00             & \textbf{48.00}    & 50.00             & 10.00   & 10.00   & 12.00              \\
23             & 50.00     & 29.00     & 50.00    & 10.00    & 10.00    & 12.00       & 50.00             & \textbf{26.33} & 50.00             & 10.00   & 10.00   & 12.00              \\
24             & 50.00     & 20.00     & 38.00    & 10.00    & 10.00     & 12.00      & 50.00             & 20.00             & \textbf{36.00}    & 10.00   & 10.00   & 12.00  \\
\hline           
\end{tabular}
}
\end{table*}\

\subsection{Load Uncertainty}
\begin{figure}[t]
	\begin{minipage}[b]{\linewidth}	
	\centering
 	\includegraphics[width=\linewidth,height=4.5cm]{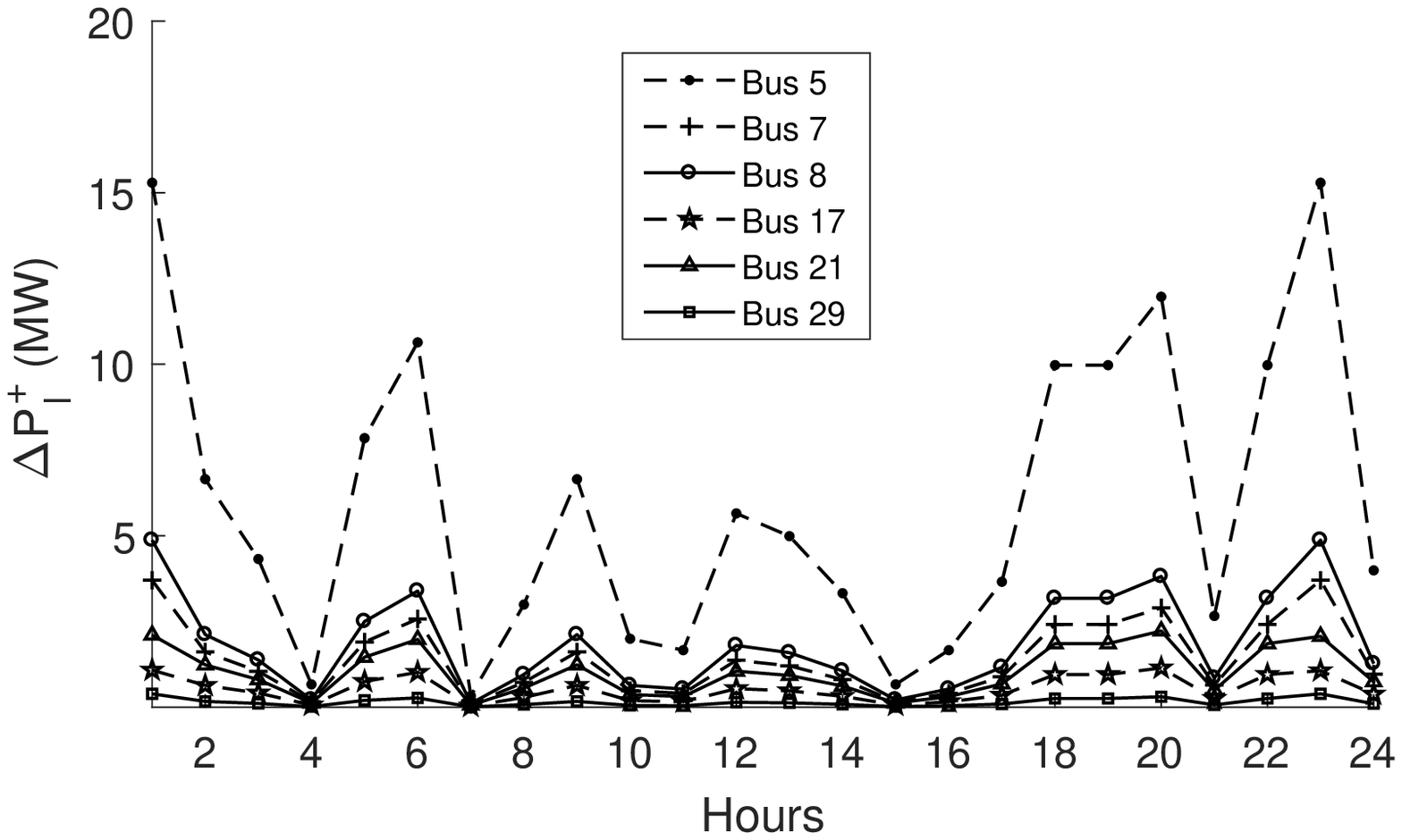}
 		\\ 
 	 	\scriptsize{(a)}
	\end{minipage}
	\begin{minipage}[b]{\linewidth}	
	\centering
 	\includegraphics[width=\linewidth,height=4.5cm]{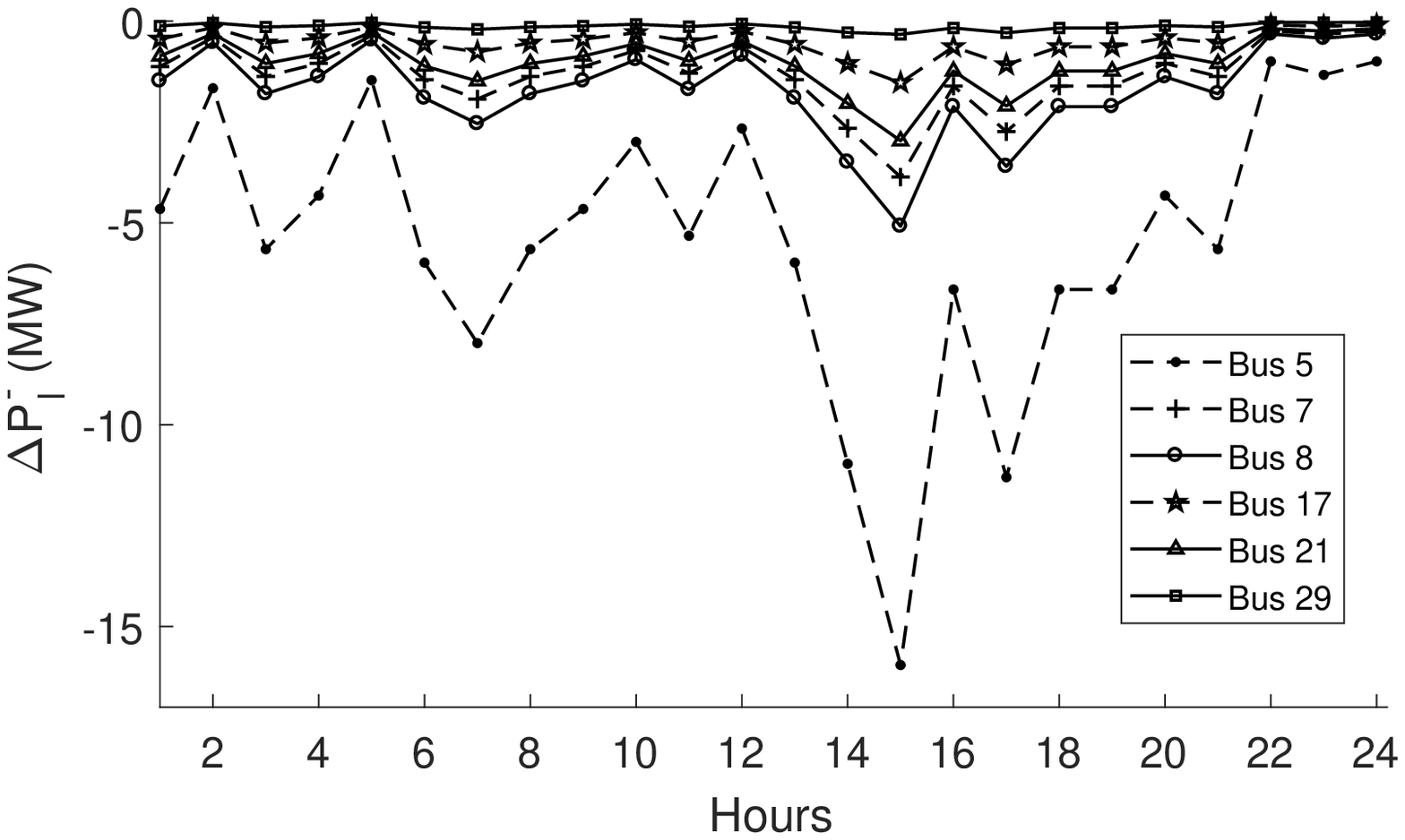}
 		\\ 
 	 	\scriptsize{(b)}
	\end{minipage}
 	\caption{Acceptable change in load (MW)}
 	\label{fig:itr_mw}
\end{figure}

The multiple simultaneous load perturbations are handled using the ITR values. These values vary in magnitude according to the present bus load and remain close to each other when seen as a fraction of the particular bus load. Hence, the normalized quantification i.e. confidence level is used in this work gives analogous results. The variation of the area under the PDF for total system load and individual bus load is quite close to each other. The acceptable load uncertainty values, in MW, are shown in Fig.\ref{fig:itr_mw}. The MW values are evidently in proportion to the load factors of different buses at which the system load is distributed among buses. The uncertainty handling capacity of the bus is dependent upon its present load.
\begin{table}[t]
  \centering
  \caption{Confidence level in load uncertainty handling}
    \begin{tabular}{c|cccccc}
    \multirow{2}[1]{*}{$\boldsymbol{\sigma}$} & \multicolumn{6}{c}{\textbf{Hours}} \\
          & \textbf{1} & \textbf{2} & \textbf{3} & \textbf{4} & \textbf{5} & \textbf{6} \\
    \hline
    \textbf{2\%} & 100.00 & 89.43 & 99.77 & 64.72 & 77.33 & 99.95 \\
    \textbf{5\%} & 95.35 & 67.18 & 80.34 & 39.50 & 57.14 & 89.74 \\
    \textbf{10\%} & 79.40 & 44.72 & 48.60 & 21.63 & 35.92 & 62.64 \\
\hline          & \multicolumn{6}{c}{\textbf{Hours}} \\
          & \textbf{7} & \textbf{8} & \textbf{9} & \textbf{10} & \textbf{11} & \textbf{12} \\
\hline    \textbf{2\%} & 59.87 & 98.42 & 99.98 & 96.83 & 95.99 & 99.38 \\
    \textbf{5\%} & 51.36 & 74.80 & 90.72 & 64.07 & 74.26 & 82.52 \\
    \textbf{10\%} & 35.39 & 45.16 & 61.74 & 35.69 & 49.68 & 54.92 \\
\hline       & \multicolumn{6}{c}{\textbf{Hours}} \\
          & \textbf{13} & \textbf{14} & \textbf{15} & \textbf{16} & \textbf{17} & \textbf{18} \\
\hline   \textbf{2\%} & 100.00 & 99.70 & 69.14 & 86.43 & 98.78 & 100.00 \\
    \textbf{5\%} & 94.46 & 86.41 & 57.93 & 62.73 & 81.30 & 94.33 \\
    \textbf{10\%} & 66.64 & 67.13 & 52.85 & 39.22 & 58.99 & 68.82 \\
\hline         & \multicolumn{6}{c}{\textbf{Hours}} \\
          & \textbf{19} & \textbf{20} & \textbf{21} & \textbf{22} & \textbf{23} & \textbf{24} \\
\hline    \textbf{2\%} & 100.00 & 99.78 & 97.72 & 78.81 & 88.49 & 85.31 \\
    \textbf{5\%} & 94.83 & 87.22 & 73.97 & 62.50 & 68.44 & 62.00 \\
    \textbf{10\%} & 69.93 & 66.09 & 45.22 & 51.41 & 59.03 & 38.83 \\
    \hline
    \end{tabular}%
  \label{tb:load_area}%
\end{table}
Therefore, absolute values cannot be used for the quantification of load uncertainty in general. The normal PDF (\ref{eq:NPDF}) is used to calculate the confidence level for uncertainty handling. The results for total system load uncertainties at three different $\sigma$ values, 2\%, 5\% and 10\% of mean are shown in the Table \ref{tb:load_area}.\

For $\sigma$=2\% the confidence is maximum at $7^{th}$ hour with 89.86\% area covered by the proposed method while there are 10 occasions having the confidence of 99\% or more. The results with 5\% standard deviation have minimum confidence level at $4^{th}$ hour while the maximum is at $1^{st}$ hour. In comparison to the previous case, values do decrease for all instances. The results with a very high $\sigma$, 10\%, is taken to illustrate the proposed method's applicability in extreme conditions. Out of all 24 instances considered, the confidence level is above 60\% at eight instances, ten instances it is between 40\% and 60\% while six instances have the confidence level below 40\% with 21.62\% (the minimum) at  $4^{th}$ hour. The $4^{th}$ hour is the only instance below 30\%. \

Further, in situations where the uncertainty of two or more PV systems or load scenarios needs to be handled, over a single base point, the SA and ITR ranges can be updated as described in the section \ref{sec:proposed}. Table \ref{tb:RUPV} shows percentage change in the values of $\Delta P_{PV}$ which can be handled further by proposed method after handling PV uncertainty at $5^{th}$ bus. The decrease in the PV output is taken as 2.25, 5, 1.5, 1.25, 4.25, 3.75, 2.5, 0.5, 1.25, and 2.75 (in MW) for $7^{th}$ to $16^{th}$ hour respectively. Results show that the decrease in $\Delta P^min_{PV}$ range is higher and occurs at 32 instances out of total 50. The increase,whenever occurred, is negligible on the other hand. $\Delta P^max_{PV}$ at all buses either increases by a large amount or remains almost same. Nevertheless, the net reduction in the uncertainty handling capacity is not more than 25\% at any bus at any instance. This shows that proposed method can be used to handle the uncertainties of multiple PV plants integrated into the system in successive fashion.\

Table \ref{tb:RULo} contains the values of the percentage change in acceptable $\Delta P_l$ (MW) after handling the same PV uncertainty. Out of all non-zero load buses, the three buses (26,29, and 30) differs at only one instance in $\Delta P^-_l$. Baring to that, all instances have the similar values for all the buses. The net change shows that the ranges can decrease and increase with comparable values. Further, the net decrease never goes beyond 24\% while net increase goes up to 28.13\%. This validates that proposed method can be used for the uncertainty handling of both PV and load at each instance in a successive manner.       
\subsection{LMP Invariability}
The variability of LMP depends upon dual variables as described in the section \ref{sec:LMP} which remains constant as the optimal basis is same. 
The difference between LMPs calculated using the MATLAB solver provided dual variables and dual variables obtained using the proposed method (\ref{eq:dual}) comes out to be in the order of $10^{-15}$. 
It proves that LMP's are not changed before and after the uncertainty in the generation/load.  

\begin{table}[t]
  \centering
  \caption{Percentage change in acceptable $\Delta P_{PV}$}
    \begin{tabular}{cc|rrrrr}
    \hline
          & \multirow{2}[2]{*}{\textbf{Bus No.}} & \multicolumn{5}{c}{\textbf{Hours}} \\
          &       & \textbf{7} & \textbf{8} & \textbf{9} & \textbf{10} & \textbf{11} \\
    \hline
    \multicolumn{1}{c}{\multirow{5}[1]{*}{\begin{sideways}{$ \Delta P_{PV}^{max}(j) $}\end{sideways}}} & \textbf{5} & 1.04  & 13.24 & 35.71 & 16.67 & 7.81 \\
          & \textbf{18} & 1.04  & 13.24 & 35.71 & 16.67 & 7.81 \\
          & \textbf{26} & 0.02  & 0.20  & 0.46  & 16.67 & 0.00 \\
          & \textbf{29} & 0.02  & 0.20  & 0.46  & 0.00  & 0.00 \\
          & \textbf{30} & 0.01  & 0.08  & 0.20  & 0.00  & 0.00 \\ \hline
    \multicolumn{1}{c}{\multirow{5}[1]{*}{\begin{sideways}{$ \Delta P_{PV}^{min}(j) $}\end{sideways}}} & \textbf{5} & -25.00 & -25.00 & -25.00 & -25.00 & -25.00 \\
          & \textbf{18} & -25.00 & -25.00 & 0.52  & -25.00 & -25.00 \\
          & \textbf{26} & -25.00 & 0.08  & 0.20  & 0.00  & 0.00 \\
          & \textbf{29} & -25.00 & -14.96 & 0.00  & -25.00 & -25.00 \\
          & \textbf{30} & -25.00 & -25.00 & -21.51 & -25.00 & -25.00 \\
    \hline
          & \multirow{2}[2]{*}{\textbf{Bus No.}} & \multicolumn{5}{c}{\textbf{Hours}} \\
          &       & \textbf{12} & \textbf{13} & \textbf{14} & \textbf{15} & \textbf{16} \\
    \hline
    \multicolumn{1}{c}{\multirow{5}[1]{*}{\begin{sideways}{$ \Delta P_{PV}^{max}(j) $}\end{sideways}}} & \textbf{5} & 53.12 & 20.83 & 7.58  & 1.04  & 6.25 \\
          & \textbf{18} & 53.12 & 20.83 & -0.06 & -0.01 & 6.25 \\
          & \textbf{26} & 36.49 & 0.37  & 0.24  & 0.04  & 0.14 \\
          & \textbf{29} & 0.20  & 0.37  & 0.24  & 0.04  & 0.14 \\
          & \textbf{30} & 0.20  & 0.16  & 0.10  & 0.02  & 0.06 \\ \hline
    \multicolumn{1}{c}{\multirow{5}[1]{*}{\begin{sideways}{$ \Delta P_{PV}^{min}(j) $}\end{sideways}}} & \textbf{5} & -25.00 & -25.00 & -25.00 & -25.00 & -25.00 \\
          & \textbf{18} & 0.49  & 0.41  & -13.16 & -25.00 & -25.00 \\
          & \textbf{26} & 0.20  & 0.16  & 0.10  & -25.00 & 0.06 \\
          & \textbf{29} & 0.00  & 0.00  & 0.00  & -25.00 & -25.00 \\
          & \textbf{30} & -25.00 & -25.00 & -25.00 & -25.00 & -25.00 \\
    \hline
    \label{tb:RUPV}%
    \end{tabular}%
\end{table}%
\begin{table}[t]
   \caption{Percentage change in acceptable $\Delta P_{l}$}
    \begin{tabular}{cc|ccccc}
    \hline
    \multirow{2}[1]{*}{\textbf{Bus No.}} & \multirow{2}[1]{*}{} & \multicolumn{5}{c}{\textbf{Hours}} \\
          &       & \textbf{7} & \textbf{8} & \textbf{9} & \textbf{10} & \textbf{11} \\
    \hline
    \multicolumn{1}{c}{\multirow{3}[1]{*}{1-24}} & $\Delta P_l^+(j)$ & -25.00 & -25.00 & -25.00 & -25.00 & -25.00 \\
          & $\Delta P^-_{l}(j)$ & 1.04  & 13.24 & 35.71 & 16.67 & 7.81 \\
          & Net   & -23.96 & -11.76 & 10.71 & -8.33 & -17.19 \\ \hline
    \multicolumn{1}{c}{\multirow{3}[1]{*}{26,29,30}} & $\Delta P^+_{l}(j)$& -25.00 & -25.00 & -25.00 & -25.00 & -25.00 \\
          & $\Delta P^-_{l}(j)$ & 1.04  & 13.24 & 35.71 & 16.67 & 7.81 \\
          & Net   & -23.96 & -11.76 & 10.71 & -8.33 & -17.19 \\
    \hline
    \multirow{2}[2]{*}{\textbf{Bus No.}} & \multirow{2}[2]{*}{} & \multicolumn{5}{c}{\textbf{Hours}} \\
          &       & \textbf{12} & \textbf{13} & \textbf{14} & \textbf{15} & \textbf{16} \\
    \hline
    \multicolumn{1}{c}{\multirow{3}[1]{*}{1-24}} & $\Delta P^+_{l}$ & -25.00 & -25.00 & -25.00 & -25.00 & -25.00 \\
          & $\Delta P^-_{l}$ & 53.13 & 20.83 & 7.58  & 1.04  & 6.25 \\
          & Net   & 28.13 & -4.17 & -17.42 & -23.96 & -18.75 \\ \hline
    \multicolumn{1}{c}{\multirow{3}[0]{*}{26,29,30}} & $\Delta P^+_{l}$ & -25.00 & -25.00 & -25.00 & -25.00 & -25.00 \\
          & $\Delta P^-_{l}$ & 53.12 & 20.83 & 7.58  & 0.00  & 6.25 \\
          & Net   & 28.12 & -4.17 & -17.42 & -25.00 & -18.75 \\ \hline
    \end{tabular}%
  \label{tb:RULo}%
\end{table}%

The time complexity of the proposed method is very less. 
It is because with the optimal solution, the optimal basis information is also available. The time taken is further reduced as the mathematical tools are concentrated upon the generation/load changes in proposed method. It takes, on an average, 1.29 ms for SA, 2.21 ms for ITR, 1.36 ms for $\beta_k$ and 3.55 ms of clock time for range update and new SA, ITR calculations. Further, the majority of time is taken during matrix inversion which can be optimized using sparse matrix operations. 

\section{Conclusion}
The proposed work presents a method for handling uncertainties in Dynamic DA-OPF using the post-optimal analysis of LP problem. The method does not require any uncertainty information for obtaining the solution and maintains optimality under uncertainties. It can handle single or multiple solar PV and load uncertainties using proposed range update method. The change in load/generation due to uncertainties is distributed among generators using proposed participation factors without re-optimization. The numerical results of confidence level show significant promise to handle uncertainties. The proposed methods can be extended for various applications and formulations like real-time economic dispatch, look-ahead optimal power flow.  
\appendices
\counterwithin{figure}{section}

\section{Uncertainty Modeling}
In the present work, the uncertainty is modeled to quantify the amount which can be handled through post-optimal analysis i.e. confidence level. Hence, the accuracy of the proposed method does not depend upon the statistical modeling of the uncertainties.  
\begin{figure}[b]
	 \begin{minipage}[b]{\linewidth}	 	
	 	\centering 	
	 \includegraphics[width=\linewidth,height=4cm]{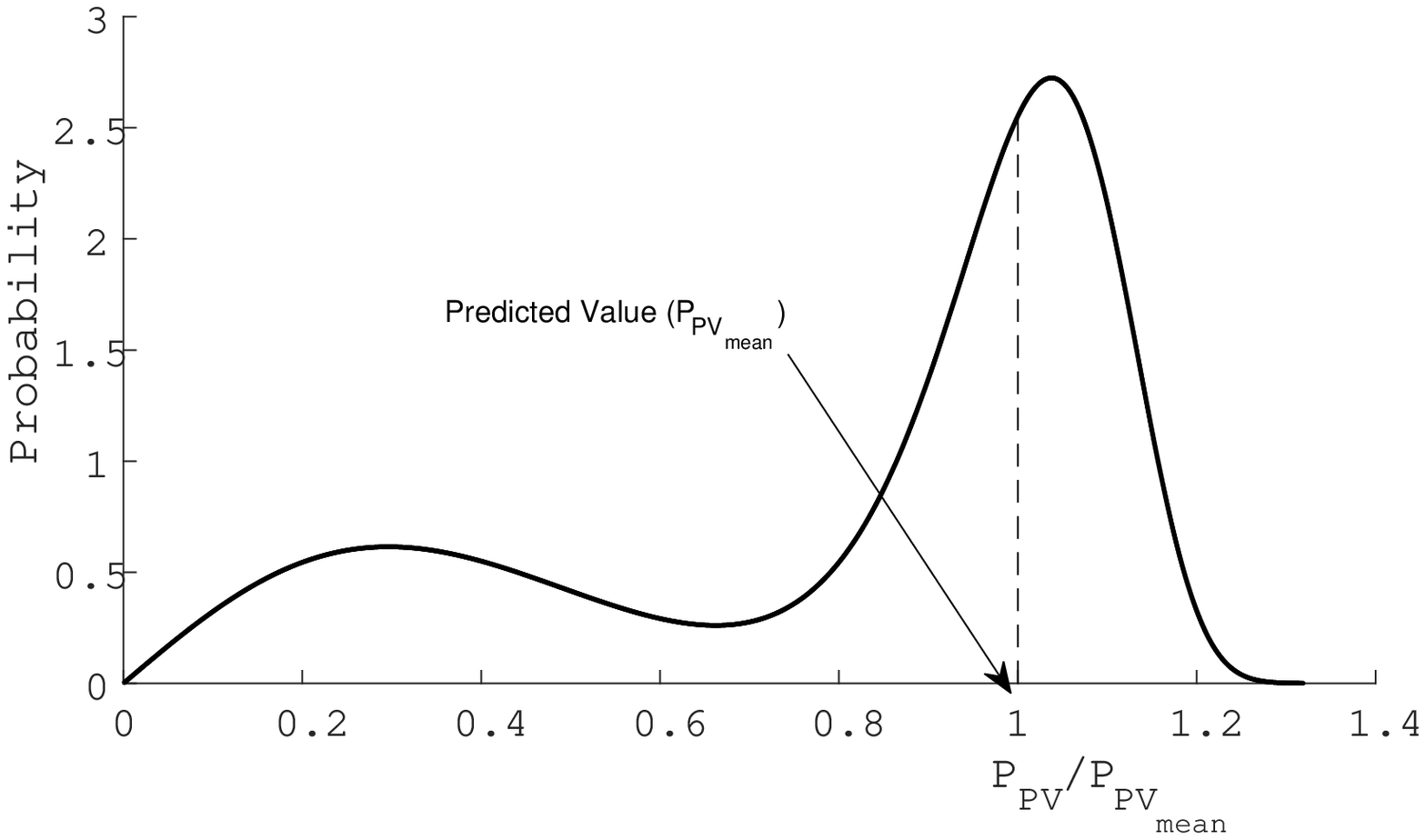}\\
	 	\scriptsize{(a) Solar PV}
     \end{minipage}
 \begin{minipage}[b]{\linewidth}
 	\centering
 	\includegraphics[width=\linewidth,height=4cm]{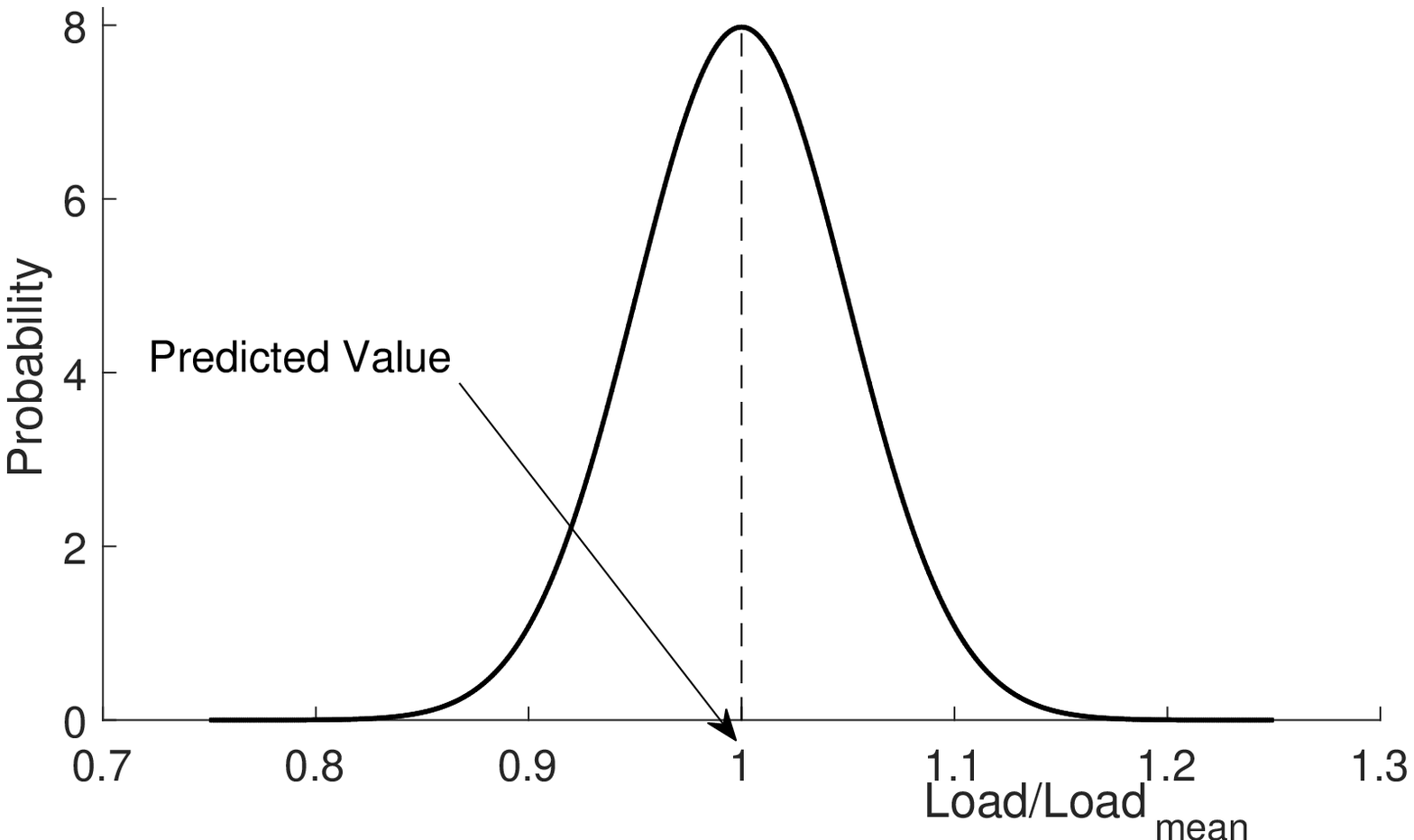}
 	\\ 
 	\scriptsize{(b) Load}
 	 \end{minipage}
\caption{An inductive distribution models }
\label{fig:un}
\end{figure}
\subsubsection{Solar Power Uncertainty}
The output of PV) plant is directly dependent upon the irradiance. The hourly distribution of irradiance, at a location, is considered as bi-modal distribution. The realization can be seen as an affine combination of two uni-modal distributions\cite{wang2011spinning}. The distribution function is modeled with Weibull PDF (with $0<g<\infty$) as:
\begin{equation}
\begin{split}
f(g)=w_1(k_1/c_1)(g/c_1)^{(k_1-1)}e^{-(g/c_1)^{k_1}}\\~+w_2(k_2/c_2)(g/c_2)^{(k_2-1)}e^{-(g/c_2)^{k_2}},
\end{split}
\end{equation}  
\ \ Where $g$ is irradiance ($kW/m^2$), $w_1$ and $w_2$ are weighted factors, $ c_1 $ and $c_2$ defines scale and $k_1$ and $k_2$ are the shape factors. Now with $\eta$ as efficiency and $S$ as capacity of the Solar PV power plant, the output power ($P_{PV}$) will be\cite{chedid1998decision}:
\begin{equation}
P_{PV}=\eta Sg
\end{equation}
\ \ By taking the prediction as mean of this distribution, the distribution of $P_{PV}$ will look as shown in Fig.\ref{fig:un}(a).

\subsubsection{Load Uncertainty}
Even though over the years load forecasting has become very accurate, the possibility of uncertainty cannot be overruled. In this paper, the Normal distribution is used to model the load uncertainty \cite{bo2009impact}. The prediction is taken as mean and three different standard deviations (2\%, 5\%, and 10\%) are used for uncertainty quantifications. The normal PDF of the load is given as:
\begin{equation}\label{eq:NPDF}
f(l)=\dfrac{1}{\sigma\sqrt{2\pi}}\times\exp\Bigg\{{\dfrac{-(l-\mu)^2}{2\sigma^2}}\Bigg\}
\end{equation}

Here, $\mu$ is the mean and $\sigma$ is standard deviation of the normal distribution. A normalized distribution with $\sigma$ as 5\% is shown in Fig.\ref{fig:un}(b). 

\ifCLASSOPTIONcaptionsoff
  \newpage
\fi




%

\bibliographystyle{IEEEtran}
\bibliography{ref}
\include{thebibliograpgy}

\end{document}